\input amstex
\documentstyle{conm-p-mod}
\magnification=1200
\NoBlackBoxes
\define\up{\phantom{$^{1^1}$}}
\define\dn{\phantom{$_{j_j}$}}
\define\hs{\hskip.7pt}
\define\nh{\hskip-.7pt}
\def\bbbc{{\mathchoice {\setbox0=\hbox{$\displaystyle\text{\rm C}$}\hbox{\hbox
to0pt{\kern0.4\wd0\vrule height0.9\ht0\hss}\box0}}
{\setbox0=\hbox{$\textstyle\text{\rm C}$}\hbox{\hbox
to0pt{\kern0.4\wd0\vrule height0.9\ht0\hss}\box0}}
{\setbox0=\hbox{$\scriptstyle\text{\rm C}$}\hbox{\hbox
to0pt{\kern0.4\wd0\vrule height0.9\ht0\hss}\box0}}
{\setbox0=\hbox{$\scriptscriptstyle\text{\rm C}$}\hbox{\hbox
to0pt{\kern0.4\wd0\vrule height0.9\ht0\hss}\box0}}}}
\font\smallit=ptmri at 10pt
\font\smallbf=ptmbo at 10pt
\topmatter
\topspace{-31pt}
\title{Geometry of the Standard Model}
\endtitle
\author Andrzej Derdzinski\endauthor
\address Department of Mathematics, The Ohio State University, 
Columbus, Ohio 43210\endaddress
\email andrzej\@{}math.ohio-state.edu\endemail
\endtopmatter
\voffset=-20pt
\hoffset=20pt
\document
\subhead 0. Introduction\endsubhead The {\it Standard Model\/} of particles 
and interactions is the currently-accepted theory of elementary particles. It 
can be naturally divided into the {\it classical part\hs}, a description of 
which is possible in the language of vector bundles over the spacetime and 
operations on them, and a {\it field quantization\/} procedure that transforms 
the classical part into a reasonable model of physical reality.\par This note 
covers only the {\it classical\/} part of the Standard Model. Similar but more 
detailed expositions of this topic can be found in the following texts:
\vskip4pt
{\eightpoint
\ref\by A. Derdzinski\book Geometry of the Standard Model of Elementary 
Particles\publ Springer-Verlag\publaddr Berlin-Heidelberg-New York\bookinfo 
Texts and Monographs in Physics\yr1992\endref
\ref\by A. Derdzinski\paper Geometry of elementary particles\jour Proceedings 
of Symposia in Pure Mathematics\vol54\yr1993\pages(edited by R. E. 
Greene and S.-T. Yau), Part 2, 157--171\endref}
\subhead 1. Interactions\endsubhead 
Aside from gravity, the known kinds of particle interactions, ordered by 
decreasing strength, are the \it strong, electromagnetic, \rm and \it weak \rm 
forces. The latter two may  be combined into the {\it electroweak 
interaction\/} (\S5).\par The {\it strength\/} of an interaction amounts to 
the probability of its occurrence in the given circumstances.
\subhead 2. Taxonomy of particle species\endsubhead(See also \S3, \S 5.)
\vskip2pt
{\eightpoint
\vbox
{\hbox{\hskip90ptknown particles}
\vskip1pt
\hbox{$\overbrace{\hskip265pt}$}
\vskip-2pt
\hbox{interaction carriers:\hskip135ptmatter particles}
\vskip1pt
\hbox{$8\,$ species of gluons,\hskip20pt
$\overbrace{\hskip245pt}$}
\vskip1pt
\hbox{$\gamma\,$ (the photon),\hskip32pt
\vbox{\vskip-3pt\hbox{leptons}\vskip3pt}\hskip110pt
\vbox{\vskip-3pt\hbox{hadrons}\vskip3pt}}
\vskip1pt
\hbox{$\text{\rm W}^-\!,\hs\text{\rm W}^+\!,\hs\text{\rm Z}^0$\hskip10pt
$\overbrace{\hskip105pt}$\hskip28pt
$\overbrace{\hskip110pt}$}
\hbox{\hskip55pt\vbox{\hbox{leptons proper,}
\hbox{$3$ generations:}
\hbox{$\text{\rm e},\nu_{\text{\rm {e}}}$}
\hbox{$\mu,\nu_\mu$}
\hbox{$\tau,\nu_\tau$}}\hskip8pt
\vbox{\hbox{antileptons,}
\hbox{$3$ generations:}
\hbox{$\text{\rm e}^+\!,\overline\nu_{\text{\rm e}}$}
\hbox{$\mu^+\!,\overline\nu_\mu$}
\hbox{$\tau^+\!,\overline\nu_\tau$}}
\hskip10pt
\vbox{\hbox{mesons}\hbox{$(>100)$}}\hskip10pt
\vbox{\hbox{\hskip40ptbaryons}\vskip3pt
\hbox{$\overbrace{\hskip106pt}$}
\hbox{\vbox{\hbox{baryons proper}\hbox{$(>50)$}}\hskip5pt
\vbox{\hbox{antibaryons}\hbox{$(>50)$}}}}}}}
\subhead 3. Definitions\endsubhead{\it 
Interaction carriers} mediate interactions, {\it matter} particles do not. 
{\it Leptons} can't interact  strongly, {\it hadrons} can. {\it Mesons}  are 
bosons, {\it baryons} are fermions (see Table 4.2). Baryons naturally form 
the disjoint classes of {\it baryons proper} and {\it antibaryons}, which 
consist of each other's antiparticles (Table 4.2). The same principle 
applies to leptons.

\newpage
\subhead 4. A physics-geometry dictionary\endsubhead\par\vskip15pt
\topcaption{Table 4.1}Particles and bundles\endcaption
\vskip3pt
{\eightpoint{
\vbox{\offinterlineskip
\hrule
\halign{&\vrule#&\strut
\quad\hfil#\hfil
\quad\cr
&\up{\sl physics}\dn&&{\sl geometry}&\cr
\noalign{\hrule}
&\up a PARTICLE species\up
&&a BUNDLE $\zeta\!$ with some geometry over the spacetime&\cr
&&&\dn$(\Cal M,g)$; the particle is {\it represented by} 
(or {\it lives in}) $\zeta$\dn&\cr
\noalign{\hrule}
&\up classical STATES of\up&&
\vbox{\hbox{SECTIONS $\,\psi\hs$ of the bundle $\,\zeta$}\vskip-5pt}&\cr
&\dn the particle\dn&&&\cr
\noalign{\hrule}
&\up EVOLUTION of the states\dn&&FIELD EQUATIONS imposed on $\psi$&\cr
\noalign{\hrule}
&\up a MATTER particle\dn&&a VECTOR bundle&\cr
\noalign{\hrule}
}}    }}
\vskip15pt
\topcaption{Table 4.2}Operations\endcaption
\vskip3pt
{\eightpoint{
\vbox{\offinterlineskip
\hrule
\halign{&\vrule#&\strut
\quad\hfil#\hfil
\quad\cr
&\up{\sl physics: operations involving}\up
&&{\sl geometry: operations on}&\cr
&\dn{\sl matter particles}\dn&&
{\sl vector bundles $\,\zeta$}&\cr
\noalign{\hrule}
&the GENERALIZATION&&\up the DIRECT SUM $\zeta_1+\ldots+\zeta_n$\up&\cr
&of $n$ given particle species&&
\dn(all $n$ species involved live here)\dn.&\cr
\noalign{\hrule}
&&&\up a natural surjective MORPHISM $\zeta_1\!\ldots\zeta_n\!\to\!\zeta$\up&
\cr
&a COMPOSITE system
&&of the TENSOR PRODUCT $\zeta_1\!\ldots\zeta_n$. It must be&\cr
&(particle)&&symmetric/skewsymmetric in any group of iden-&\cr
&\phantom{\sl physics: operations involving ll}&&\dn tical particles, then 
called {\it bosons}/{\it fermions}\dn&\cr
\noalign{\hrule}
&ANTIPARTICLE formation&&\up complex CONJUGATE $\zeta^{\!\!\!^-}$\dn&\cr
\noalign{\hrule}
}}     }}
\vskip15pt
\topcaption{Table 4.3}Interactions and gauge fields\endcaption
\vskip3pt
{\eightpoint{
\vbox{\offinterlineskip
\hrule
\halign{&\vrule#&\strut
\quad\hfil#\hfil
\quad\cr
&\up{\sl physics: interactions}\dn&&{\sl geometry: Yang-Mills fields}&\cr
\noalign{\hrule}
&&&\up a NATURAL VECTOR BUNDLE $\eta\to\Cal M$ of\up&\cr
&a FREE matter particle&&first order (the {\it free-particle bundle} of the 
given spe-&\cr
&&&\dn cies); naturality amounts to {\it direct observability}\dn&\cr
\noalign{\hrule}
&an INTERACTION&&\up a NON-NATURAL vector bundle $\delta\!\to\!\Cal M$ (the 
{\it interac-}&\cr
&\dn of some given kind\dn &&{\it tion bundle}) with some geometry,
 mainly a $G$-structure&\cr
\noalign{\hrule}
&\up CARRIERS of\up&&live in the AFFINE BUNDLE $\Cal C(\delta)$ whose sec-&\cr
&the interaction&&\dn tions are the compatible connections in $\delta$\dn&\cr
\noalign{\hrule}
&\up an INTERACT-\up&&the INTERACTING-PARTICLE BUNDLE 
$\alpha=\alpha(\delta,\eta)$,&\cr
&ING matter&&functorial in both $\delta,\eta$ and ``homogeneous linear'' in&\cr
&\dn particle\dn&&the free-particle bundle $\eta$ 
(basic example: $\alpha=\delta\eta$)&\cr
\noalign{\hrule}
}}     }}\par
\vskip15pt
Usually, neither $\alpha$ nor $\Cal C(\delta)$ is natural. This contradicts 
the obvious requirement that carriers of interactions and interacting matter 
should be directly observable. One resolves this problem by ``restoring'' 
naturality of the bundles in question using {\it bound states,} or {\it 
symmetry breaking}, as described below. (Notations: $N$ is the fibre dimension
 of the fixed interaction bundle $\delta$; an integer $k>0$ represents the 
product vector bundle $\Cal M\times\bbbc^{\hskip1ptk}\hskip-2pt$.)

\newpage
\topcaption{Table 4.4}Bound states and symmetry breaking\endcaption
\vskip4pt
{\eightpoint{
\vbox{\offinterlineskip
\hrule
\halign{&\vrule#&\strut
\quad\hfil#\hfil
\quad\cr
&\up{\sl physics}\up&&\dn{\sl geometry}\dn&\cr
\noalign{\hrule}
&\up BOUND STATES\up&&MORPHISMS of $\alpha_{_1}\!\ldots\alpha_{_n}$ 
onto NATURAL bundles,&\cr
&\dn of $n$ particles\dn&&obtained by naturally ``canceling'' 
the $\delta$-related factors&\cr
\noalign{\hrule}
&\vbox{\hbox{SYMMETRY BREAKING}\vskip-5.5pt}&&\up selection of an ADDITIONAL 
STRUCTURE\up&\cr
&&&\dn in $\delta$, leading to reduction of $G$ to a subgroup\dn&\cr
\noalign{\hrule}
&&&\up TRIVIALIZATION of $\delta$, so $\delta=N$ and, e.g., 
$\delta\eta=N\hskip-1.3pt\eta$\up
&\cr
&FORMAL symmetry breaking&&(the interacting particle comes in $N$
separate versions),&\cr
&(a thought experiment)&&while $\Cal C(\delta)=(\dim G)T^*$, i.e., 
the carriers appear as&\cr
&&&\dn$\dim G$ species of matter particles living in 
$T^*\!=T^*\hskip-2.4pt\Cal M$\dn&\cr
\noalign{\hrule}
&\up SPONTANEOUS symmetry\up&&&\cr
&breaking (in nature, for in-&&Example: the ELECTROWEAK MODEL (\S $5$).&\cr
&\dn teractions of low strength)\dn&&&\cr
\noalign{\hrule}
}}     }}
\vskip20pt
\subhead 5. The standard model\endsubhead\par
\vskip12pt
\topcaption{Table 5.1}Geometry of interactions\endcaption
\vskip4pt
{\eightpoint{
\vbox{\offinterlineskip
\hrule
\halign{&\vrule#&\strut
\quad\hfil#\hfil
\quad\cr
&{\sl inter-}&&\up ELECTRO-\up&&\vbox{\hbox{ELECTROWEAK}\vskip-5.5pt}&
&\vbox{\hbox{STRONG}\vskip-5.5pt}&\cr
&{\sl action}&&\dn MAGNETIC\dn &&&&&\cr
\noalign{\hrule}
&\vbox{\hbox{\sl credits}\vskip-5pt}&&\up Weyl,\up&&Glashow, Salam, Weinberg,&
&Gell-Mann, Zweig,&\cr
&&&\dn1929\dn&&1961--1967&&1964&\cr
\noalign{\hrule}
&{\sl }&&
\up The possibility of\up&&
The model describes one gen-&&
&\cr
&{\sl }&&
a unified descrip-&&
eration of (anti)leptons (\S $2$)&&
Hadrons appear&\cr
&{\sl }&&
tion of electromag-&&
at a time. Choose, e.g.,  $\text{\rm {e}},\nu_{\text{\rm {e}}}\,$:&&
as composites&\cr
&{\sl }&&
netism for all par-&&
the electron and electronic neu-&&
of {\it quarks} and&\cr
&{\sl com-}&&
ticles expresses&&
trino. Their free-particle bun-&&
{\it antiquarks}&\cr
&{\sl ments}&&
the fact that the&&
dles are: $\sigma$ for {\rm e} and $\sigma\!_{_{\text{L}}}$ 
for $\nu_{\text{e}}$,&&
(abbreviation:&\cr
&{\sl }&&
electric charge is&&
where $\sigma$ denotes a fixed Dirac&&
q, \B q), coming&\cr
&{\sl }&&
{\it quantized}, i.e., oc-&&
spinor bundle, $\Cal M$ is assum-&&
in several {\it fla-}&\cr
&{\sl }&&
curs in multiples&&
ed orientable, and $\sigma\!=\!\sigma\!_{_{\text{L}}}\hskip-3pt+\!\sigma\!_
{_{\text{R}}}$&&
{\it vors} (species).&\cr
&{\sl }&&
\dn of a fixed amount.\dn&&
(Weyl spinor bundles),  $\sigma\!_{_{\text{R}}}\!=
\overline{\sigma\!_{_{\text{L}}}}.$&&
&\cr
\noalign{\hrule}
&\up$G,\, \delta$\dn&&$G\!=\!\text{{\rm U}}(1),\,\delta\!=\!\lambda$&&
$G\!=\!\text{{\rm U}}(2),\,\delta\!=\!\iota$&&
$G\!=\!\text{{\rm SU}}(3),\,\delta\!=\!\rho$&\cr 
\noalign{\hrule}
&\vbox{\hbox{\sl what $\delta$ is}\vskip-5pt}&&\up a complex\up&&a complex&&a 
complex&\cr
&{\sl }&&\dn line bundle\dn&&plane bundle&&$3$-space bundle&\cr
\noalign{\hrule}
&{\sl geometry}
&&\up a Hermitian\up&&
a Hermitian&&
$\langle\,,\rangle$ and a unit sec-&\cr
&\dn{\sl of $\delta$\dn}
&&fibre metric $\langle\,,\rangle$&&fibre metric $\langle\,,\rangle$&&
tion $\Theta$ of $[\rho^*]^{\wedge3}$&\cr
\noalign{\hrule}
&{\sl free-par-}&&&&\up a fixed Dirac spinor\up&&
$\sigma$ for quarks,&\cr
&{\sl ticle}&&any $\eta$&&bundle $\sigma$ for the whole
&&$\overline{\sigma}$ for&\cr
&{\sl bundle}&&&&
\dn generation $\text{\rm {e}},\nu_{\text{\rm {e}}}$\dn&&
antiquarks&\cr
\noalign{\hrule}
&{\sl inter-}&&$\alpha\!=\!\lambda^k\eta\,$ if parti-&
&\up$\alpha=\iota\sigma\!_{_{\text{L}}}+\iota^{\wedge2}\sigma\!_{_{\text{R}}}$
\up&
&$\alpha=\rho\sigma$&\cr
&{\sl acting}&&cle carries {\it k} units&&or, if neutrinos are massive,&
&(for quarks)&\cr
&{\sl particle}&&of electron charge&&even simpler:&
&$\alpha=\overline{\rho}\overline{\sigma}$&\cr
&{\sl bundle}&&(with $\lambda^{-k}=\overline{\lambda^k}$)\dn&&
$\alpha=\iota\sigma$&&(antiquarks)&\cr
\noalign{\hrule}
}}     }}\rm\par

\newpage
\topcaption{Table 5.2}Bound states and symmetry breaking in the standard 
model\endcaption
\vskip7pt
{\eightpoint{
\vbox{\offinterlineskip
\hrule
\halign{&\vrule#&\strut
\quad\hfil#\hfil
\quad\cr
&\up{\sl inter-}\up&&ELECTRO-&&\vbox{\hbox{ELECTROWEAK}\vskip-5.5pt}&
&\vbox{\hbox{STRONG}\vskip-5.5pt}&\cr
&\dn{\sl action}\dn&&MAGNETIC&&&&&\cr
\noalign{\hrule}
&{\sl bound}&&\up only if\up&&&&cancel $\rho$ factors&\cr
&{\sl states:}&&
$\sum_{j=1}^n\!\text{\it k}_{_j}=0$&&&&
by $\langle\,,\rangle\!:\!\rho\overline{\rho}\!\to\!1$,&\cr
&{\sl $\alpha_{_1}\!\ldots\alpha_{_n}$}&&
(electrically neu-&&&&$\Theta\!:\!\rho^3\!\to\!1$, or&\cr
&{\sl $\downarrow$}&&
tral systems,&&\vbox{\hbox{none of interest}\vskip-5pt}
&&$\overline{\Theta}\!:\!\overline{\rho}^3\!\!\to\!\!1$, getting&\cr
&{\sl $\zeta$}&&
e.g., atoms),&&&&q\B q pairs (mesons),&\cr
&{\sl (where both}&&
as $\lambda\overline{\lambda}\!=\!1$ and&&&&q triples (bary-&\cr
&{\sl $\zeta$ and $\downarrow$}&&
$\lambda^{k_1}\!\!\ldots\lambda^{k_n}\!=$&&&&ons proper), \B q tri-&\cr
&{\sl are natural\hs)}&&$\lambda^{\sum_j\!k_j}\!$ under 
$\langle\,,\rangle$&&&&ples\dn\hskip3pt(antibaryons)&\cr
\noalign{\hrule}
&{\sl }&&
\up$\lambda\!=\!1,\alpha\!=\!\eta,$\up&&&
&$\rho\!=\!3,\,\,\alpha\!=\!3\sigma\,$ or&\cr
&{\sl formal}&&$\Cal C(\lambda)\!=\!T^*\!.$ Mat-&&&&
$3\overline{\sigma}\text{:}\!$ each q,\hs\B q flavor&\cr
&{\sl symmetry}&&ter: same as free.&&
\vbox{\hbox{of no interest}\vskip-5pt}&&comes in $3\!$ {\it colors}.&\cr
&{\sl breaking}&&Carriers: just&&&&
As $\Cal C(\rho)\!=\!8T^*\!\!,\!$ the&\cr
&{\sl }&&one species,&&&&carriers appear as&\cr
&{\sl }&&\dn the {\it photon} $\gamma$.\dn&&&&
$8$ species of {\it gluons.}&\cr
\noalign{\hrule}
&{\sl }&&&&\up Choice of a section $\phi$ of $\iota$&&&\cr
&{\sl }&&&&with $\!|\phi|\!=\!\text{\rm constant}\!>\!0\!$ re-&&&\cr
&{\sl }&&&&duces $\!\text{\rm U}(2)\!$ to $\!\text{\rm U}(1).\!$ Call 
$\lambda\!$&&&\cr
&{\sl }&&&&$=\!\phi^{\perp}\!$ the {\it electromagnetic-}&&&\cr
&{\sl }&&&&{\it interaction bundle,} so $\iota\!=\!$&&&\cr
&{\sl }&&&&$1\!+\!\lambda,1\!=\!\text{Span}\,\phi,
\iota^{\wedge2}\!=\!\lambda$.&&&\cr
&{\sl spontaneous}&&none:&&
Thus, $\alpha\!=\!\sigma_{_{\text{L}}}\!+\!\lambda\sigma$ describ-&&none:&\cr
&{\sl symmetry}&&too strong&&es $\text{\rm e},\nu_{\text{\rm e}}$ 
with their correct&&much&\cr
&{\sl breaking}&&&&charges, and the summands&&too strong&\cr
&&&&&of $\Cal C(\iota)\!=\!\Cal C(\lambda)\!+\nh\lambda T^*\!\!+\nh T^*$&&&\cr
&&&&&represent the carriers: the&&&\cr
&&&&&{\it photon} $\gamma$, and the massive,&&&\cr
&&&&&matterlike {\it weak-interac-}&&&\cr
&&&&&{\it tion carriers,} $\!\text{\rm W}^{\pm}\!$ (charg-&&&\cr
&&&&&ed), and $\text{\rm Z}$ (neutral).\dn&&&\cr
\noalign{\hrule}
}}    }}\rm
\vskip25pt
\subhead 6. Coupling constants and the Weinberg angle\endsubhead
An additional ingredient of the geometry of any interaction bundle $\delta$ 
is provided by a fixed {\it natural fibre metric} $(,)$ in $\Cal C(\delta)$, 
obtained in the obvious way from a biinvariant metric on $G$. Since $G=\text
{U}(2)$ is reducible, the freedom in choosing $(,)$ for the electroweak model
 involves not merely a scale factor (referred to as a {\it coupling 
constant}), but also an angular parameter (the {\it Weinberg angle}).\par
 The latter is physically meaningful, since the decomposition of $\Cal C
(\iota)$ in Table 5.2 is $(,)$-orthogonal. In general, the coupling constant
 of $(,)$ also has a physical interpretation, namely in terms of the {\it 
strength} of the interaction described by $\delta$.

\enddocument